\documentclass[12pt,leqno]{article}
\usepackage{amsmath, amsthm, amsfonts, amssymb, color}
\usepackage{mathrsfs}
\usepackage[notcite,notref]{showkeys}

\newtheorem{thm}{Theorem}[section]

\newtheorem{lem}[thm]{Lemma}

\newtheorem{Remark}[thm]{Remark}
\theoremstyle{definition}

\newcommand{\scr}[1]{\mathscr #1}
\definecolor{wco}{rgb}{0.5,0.2,0.3}

\numberwithin{equation}{section}

\newcommand{\ua}{\uparrow}

\title{{\bf Characterizations of the upper bound of Bakry-Emery curvature}\footnote{Supported in
 part by  NNSFC (11371099).} }
\author{
{\bf   Bo Wu}\\
\footnotesize { School  of Mathematical Sciences, Fudan
University, Shanghai 200433, China}\\
 \footnotesize{Institute for Applied Mathematics, University of Bonn,
Endenicher Allee 60, 53115 Bonn, Germany}\\
\footnotesize{wubo@fudan.edu.cn, wu@iam.uni-bonn.de}}

\date{}
\begin{document}
\maketitle

\def\R{\mathbb R} \def\EE{\mathbb E} \def\Z{\mathbb Z} \def\ff{\frac} \def\ss{\sqrt}
\def\H{\mathbb H}
\def\dd{\delta} \def\DD{\Delta} \def\vv{\varepsilon} \def\rr{\rho}
\def\<{\langle} \def\>{\rangle} \def\GG{\Gamma} \def\gg{\gamma}
\def\ll{\lambda} \def\LL{\Lambda} \def\nn{\nabla} \def\pp{\partial}
\def\d{\text{\rm{d}}} \def\Id{\text{\rm{Id}}}\def\loc{\text{\rm{loc}}} \def\bb{\beta} \def\aa{\alpha} \def\D{\scr D}
\def\E{\scr E} \def\si{\sigma} \def\ess{\text{\rm{ess}}}
\def\beg{\begin} \def\beq{\beg}  \def\F{\scr F}
\def\Ric{\text{\rm{Ric}}}
\def\Var{\text{\rm{Var}}}
\def\Ent{\text{\rm{Ent}}}
\def\Hess{\text{\rm{Hess}}}\def\B{\scr B}
\def\e{\text{\rm{e}}} \def\ua{\underline a} \def\OO{\Omega} \def\b{\mathbf b}
\def\oo{\omega}     \def\tt{\tilde} \def\Ric{\text{\rm{Ric}}}
\def\cut{\text{\rm{cut}}} \def\P{\mathbb P} \def\ifn{I_n(f^{\bigotimes n})}
\def\fff{f(x_1)\dots f(x_n)} \def\ifm{I_m(g^{\bigotimes m})} \def\ee{\varepsilon}
\def\C{\scr C}
\def\M{\scr M}\def\ll{\lambda}
\def\X{\scr X}
\def\T{\scr T}
\def\A{\mathbf A}
\def\LL{\scr L}\def\LLL{\Lambda}
\def\gap{\mathbf{gap}}
\def\div{\text{\rm div}}
\def\Lip{\text{\rm Lip}}
\def\dist{\text{\rm dist}}
\def\cut{\text{\rm cut}}
\def\supp{\text{\rm supp}}
\def\Cov{\text{\rm Cov}}
\def\Dom{\text{\rm Dom}}
\def\Cap{\text{\rm Cap}}\def\II{{\mathbb I}}\def\beq{\beg{equation}}
\def\sect{\text{\rm sect}}\def\H{\mathbb H}

\begin{abstract}
In this paper, we will present some characterizations for the upper bound of the Bakry-Emery curvature on a Riemannian manifold by using functional inequalities on path space. Moreover, some characterizations for general lower and upper bounds of Ricci curvature are also given, which extends the recent results derived by Naber \cite{N} and Wang-Wu\cite{WW}.  A crucial point of the present study is to use the symmetrical idea for the lower and upper bounds of Ricci curvature, and the localization technique.
\end{abstract}

\vskip 2cm

\section{Introduction}\label{sect1}
Let $(M,g)$ be an $n$-dimension complete Riemannian manifold and $Z$ be a $C^1$-vector field on $M$. Consider the Bakry-Emery curvature $\Ric_Z:=\Ric+\nabla Z$ for the Witten Laplace $L:=\Delta_M-Z$, where $\Delta_M$ is the Laplace operator on $M$. The bounded property of $\Ric_Z$ play a very important role on the field of analysis. Specially, there are a number of equivalent characterizations for the lower
bound of $\Ric_Z$ by functional inequalities of the (Neumann) semigroup
generated by $L$ in \cite{W1}. However, for the associated upper bound, there are not almost any results. In this paper, there are two goals: one is to present some characterizations for general lower and upper bounds of $\Ric_Z$ by using the symmetrical idea; Based on the above result and the localization technique, the other one is that we will give some equivalent characterizations for the upper bound of $\Ric_Z$. The first result extends the recent results derived by Naber \cite{N} and Wang-Wu\cite{WW}. The motivation of this work comes from the following observation: The upper and lower bounds of the Bakry-Emery curvature $\Ric_Z$ at any point $x\in M$ are in essence actually determined by the distribution properties of all paths near this point.

Before moving on, let us recall some notation on path space. For any $T>0$ and $x\in M$, the based path space and the path space over $M$ are defined respectively by
$$W^T(M):=C([0,T]: M)$$
and
$$W_x^T(M):=\big\{\gamma\in W^T(M): \gamma_0=x\big\}.$$
Let $d$ be the distance function on $M$.  Then $W^T(M)$ is a Polish
space under the uniform distance
$$\rho(\gamma,\sigma):=\displaystyle\sup_{t\in[0,T]}d(\gamma_t,\sigma_t),\quad\gamma,\sigma\in W^T(M).$$
In particular, let $\rho_x(\gamma):=\rho(\gamma,x), \gamma\in W^T(M)$ be the distance function on $W^T(M)$ starting from some  fixed $x\in M$.

Denote by $O_x(M)$ be the orthonormal frame bundle at $x\in M$, then $O(M):=\sup_{x\in M}O_x(M)$ is the orthonormal frame bundle over $M$. Let $U_t^x$ be the horizontal diffusion process on $O(M)$ associated to the horizontal lift of $L$; that is, $U_t^x$ solves
the following stochastic differential equation on $O(M)$,
\begin{equation}\label{eq1.1}
\d U_t^x=\sqrt{2}\displaystyle\sum^d_{i=1}H_i(U_t^x)\circ\d W_t^i+H_{Z}(U_t^x)\d t, \quad U_0\in O_x(M),\ \ t \in [0,\zeta),
\end{equation}
where $W_t=(W_t^1,\cdots,W_t^d)$ is the $d$-dimensional Brownian motion, $\{H_i\}_{i=1}^n: TM\rightarrow TO(M)$ is a standard orthonormal basis
of horizontal vector fields on $O(M)$, $H_Z$ is the horizontal vector fields of $Z$ and $\zeta$ is the maximal time of the solution. Let $\pi: O(M) \rightarrow M$ be the canonical projection.
Then $X^x_t:=\pi(U^x_t),t<\zeta$ is the $L$-diffusion process with initial point $x$.  Throughout this paper, besides the completeness of $M$,
we assume further that $X^x_t$ is non-explosive, i.e., $\zeta=\infty$, a.s..  Let $P_t$ be the semigroup generated by $L$, that is
$$P_tf(x)= \EE f(X_t^x),\ \ x\in M, f\in \B_b(M), t\ge 0.$$

For any $T>0$,
define the Cameron-Martin space
$$\mathbb{H}=\left\{h\in C([0,T];\mathbb{R}^d): h(0)=0,
\|h\|^2_{\mathbb{H}}:=\int_0^T|h_s'|^2\d s<\infty\right\},$$ which is
a separable Hilbert space under $\<h,g\>_\H:=\int_0^T \<h_s',g_s'\>\d s,\
h,g\in \mathbb H.$ Here $\F C_T^\infty$ denotes the class of cylindric smooth functions on the path space $W^T(M)$ defined as
\begin{equation}\label{eq1.2}\aligned \F C_T^\infty:=\Big\{F(\gamma)&=f(\gamma(t_1),\cdots,\gamma(t_N)):\ N\geq1,~\gamma\in W^T(M),\\&
~~~~~~~~~~~~~~~~~~~
0<t_1<t_2\cdots<t_N\leq T,f\in C_0^{Lip}(M^N)\Big\}.\endaligned\end{equation}
For each $F\in\F C_T^\infty$ with the form
$F(\gamma):=f\big(\gamma(t_1),\cdots,\gamma(t_N)\big)$ and any
$h\in\mathbb{H}$, the Malliavin derivative $D_h F$ is given by
\begin{equation}\label{eq1.3}
D_hF(X^x_{[0,T]})=\displaystyle\sum^N_{i=1}\left\<\nabla_i
f\big(X^x_{t_1},\cdots,X^x_{t_N}\big),U^x_{t_i}h(t_i)\right\>_{T_{X^x_{t_i}}M},
\end{equation}
where $\nabla_i$ is the (distributional) gradient operator for the $i$-th
component on $M^N$. The Riesz representation theorem implies that there exists a gradient function $DF\in \H$ such that
$$\big\<DF(X^x_{[0,T]}),h\big\>_{\H}=D_hF(X^x_{[0,T]}),\quad h\in \mathbb{H}.$$
In particular, if $F$ has the above form, we have
\begin{equation}\label{eq1.4}
\dot{D}_sF(X^x_{[0,T]})(s):=\frac{\d \big [DF(X^x_{[0,T]})(s)\big]}{\d s}=\sum_{t_i>s}(U^x_{t_i})^{-1}\nabla_if\big(X^x_{t_1},\cdots,X^x_{t_N}\big).
\end{equation}

On behalf of stating our main result, we also introduce some additional notation. For any $x\in M$ and $R>0$, the closed ball of radius $R$ around $x$ is denoted by $B_R(x)$, and the local cylindric functions is defined by
\begin{equation}\label{eq1.5}\aligned
\F C^{x,R,T}_{loc}&=\Big\{f\left(\gamma(t_1),\cdots,\gamma(t_N)\right)l\left(\rho_x\left(\gamma\right)\right): \\&
0<t_1<\cdots <t_N<T, \ f\in C^{Lip}_0(M^N), l\in C_o^\infty(\R),\supp(l)\subset (-\infty,R)\Big\}.
\endaligned\end{equation}
Let $K_1,K_2$ be any two continuous functions on $M$ with $K_1\geq K_2$, we introduce the following random measure on $[0,T]$:
\beq\label{MU}
  \mu^{K_1,K_2}_{x,T}(\d s):=\e^{\frac{\int^s_0 K_1(X_u^x)-K_2(X_u^x)\d u}{2}}\frac{K_1(X_s^x)-K_2(X_s^x)}{2}\d s \end{equation}
and the measurable function on $W^T(M)$:
\begin{equation}\label{eq1.7}A_{s,t}^{K_1,K_2}(X_{[0,T]}^x):=\e^{-\frac{\int^t_s K_1(X_u^x)+K_2(X_u^x)\d u}{2}}.\end{equation}
For simplicity, let $A_{t}^{K_1,K_2}(X_{[0,T]}^x):=A_{0,t}^{K_1,K_2}(X_{[0,T]}^x)$. Similar to the previous argument, the Riesz representation theorem implies that for every Lipschitz function $F$ on $W^T(M)$, there exists an unique gradient $D^{K_1,K_2}F$ such that
\begin{equation*}D_hF(X_{[0,T]}^x)=\left\langle D^{K_1,K_2}F(X_{[0,T]}^x), A_\cdot^{K_1,K_2}(X_{[0,T]}^x)h\right\rangle_{\H}, \quad \P-\text{a.s.}.\end{equation*}
In particular, if $F\in \F C_T^\infty$ with the form $F(\gamma)=f(\gamma(t_1),\cdots,\gamma(t_N))$, then
\beq\label{MG}\beg{split} &\dot{D}^{K_1,K_2}_sF(X_{[0,T]}^x)= \sum_{t_i> s}A_{t_i}^{K_1,K_2}(X_{[0,T]}^x)(U_{t_i}^x)^{-1}\nn_if\big(X_{t_1}^x,
\cdots,X_{t_N}^x\big),\ \ s\in [0,T].\end{split} \end{equation}
and $D^{K_1,K_2}F\in \H$.  In fact,
$$\aligned\left\langle D^{K_1,K_2}F, h\right\rangle_{\H}
&=\int^T_0\langle\dot{D}^{K_1,K_2}_sF,h'_s\rangle \d s\\&=\sum_{i=1}^NA_{t_i}^{K_1,K_2}
\left\langle (U_{t_i}^x)^{-1}\nn_if,h_{t_i}\right\rangle
=\sum_{i=1}^N\left\langle (U_{t_i}^x)^{-1}\nn_if,A_{t_i}^{K_1,K_2}h_{t_i}\right\rangle\\
&=\sum_{i=1}^N\int^{t_i}_0\left\langle
(U_{t_i}^x)^{-1}\nn_if,\big(A_s^{K_1,K_2} h_s\big)'\right\rangle \d s\\&=\int^T_0\left\langle
\sum_{i=1}^N1_{\{s\leq t_i\}}(U_{t_i}^x)^{-1}\nn_if,\big(A_s^{K_1,K_2}\big)'\right\rangle \d s\\
&=\int^T_0\left\langle\dot{D}_sF,\Big(A_{s}^{K_1,K_2} h_s\Big)'\right\rangle \d s=\left\<DF, A_\cdot^{K_1,K_2}h\right\>_{\H}.
\endaligned$$
Specially, when $K_2=-K_1$, we have $D^{K_1,K_2}F=DF$. Actually, $D^{K_1,K_2}F$ may be interpreted as a new gradient with $(K_1,K_2)$-weight, we will apply it to characterize the non-symmetric upper and lower bounds of the Bakry-Emery curvature.

For any $t\in [0,T]$, define respectively the gradient
\begin{equation}\label{eq1.9}\aligned&\dot{D}^{K_1,K_2}_{t,s}F(X_{[0,T]}^x)= \sum_{t_i> s}A_{t,t_i}^{K_1,K_2}(X_{[0,T]}^x)(U_{t_i}^x)^{-1}\nn_if\big(X_{t_1}^x,
\cdots,X_{t_N}^x\big), \quad s\geq t.\endaligned\end{equation}
and the energy form
$$\aligned\E_{t,T}^{K_1,K_2}(F,F)
 &= \EE\bigg\{\Big(1+\mu^{K_1,K_2}_{x,T}([t,T])\Big)
  \bigg(\big|\dot{D}^{K_1,K_2}_{t,t}F(X_{[0,T]}^x)\big|^2\\
 &~~~~~~~~~~~~~~~~~~~~~~+\int^{T}_{t}\big|\dot{D}^{K_1,K_2}_{t,s}F(X_{[0,T]}^x)\big|^2\mu^{K_1,K_2}_{x,T}(\d s)\bigg)\bigg\}\endaligned$$ for each function $F\in \F C_T^\infty$ on $W(M)$.
For any $x\in M$ and $ R>0$, let
$$C^x_R:=\inf\left\{\Ric_Z(X,X):~X\in T_yM, ~|X|=1,~y\in B_R(x)\right\}.$$

Our main results are as below:

\beg{thm}\label{T1.1} Let $K$ be a continuous function on $M$. The following statements are equivalent each other:
\beg{enumerate}\item[$(1)$] For any $x\in M$,
\beg{align*} \Ric_Z(x)\le K(x)g_x.\end{align*}
\item[$(2)$] For any $x\in M$ and for any constants $R>0, C\leq C_R^x$(or there exist constants $R>0$ and $C\leq C_R^x$), and for each $T>0$,
\beg{align*}&\big|\nn_x \EE F(X_{[0,T]}^x)\big|\le  \EE \bigg\{\big|\dot{D}^{K,C}_0F(X_{[0,T]}^x)\big| +\int_0^T
  \big|\dot{D}^{K,C}_sF(X_{[0,T]}^x)\big| \mu^{K,C}_{x,T}(\d s)\bigg\},~~ F\in \F C^{x,R,T}_{loc}.\end{align*}
  \item[$(3)$]  For any $x\in M$ and for any constants $R>0, C\leq C_R^x$(or there exist constants $R>0$ and $C\leq C_R^x$), and for each $T>0$,
\beg{align*}& \big|\nn_x \EE F(X_{[0,T]}^x)\big|^2\le  \EE \bigg\{\big(1+ \mu^{K,C}_{x,T}([0,T])\big)\\
&~~~~~~~~~~~~~~~~~~\times\bigg( \big|\dot{D}^{K,C}_0F(X_{[0,T]}^x)\big|^2 +\int_0^T
  \big|\dot{D}^{K,C}_sF(X_{[0,T]}^x)\big|^2 \mu^{K,C}_{x,T}(\d s)\bigg)\bigg\},\quad F\in \F C^{x,R,T}_{loc}.\end{align*}
\item[$(4)$] For any $x\in M$ and for any constants $R>0, C\leq C_R^x$(or there exist constants $R>0$ and $C\leq C_R^x$), and for each $T>0$ and for any  $t_0,t_1\in [0,T]$ with $t_1>t_0$,  the following log-Sobolev inequality holds: $$\aligned
    &\EE\left[\EE\big(F^2(X_{[0,T]}^x)|\F_{t_1}\big) \log \EE(F^2(X_{[0,T]}^x)|\F_{t_1})\right]\\&-\EE\left[\EE
\big(F^2(X_{[0,T]}^x)|\F_{t_0}\big) \log \EE(F^2(X_{[0,T]}^x)|\F_{t_0})\right]
 \le 4 \int_{t_0}^{t_1} \E_{s,T}^{K,C}(F,F)\d s,~ F\in \F C^{x,R,T}_{loc}.\endaligned$$
 \item[$(5)$] For any $x\in M$ and for any constants $R>0, C\leq C_R^x$(or there exist constants $R>0$ and $C\leq C_R^x$), and for each $T>0$ and for any  $t_0,t_1\in [0,T]$ with $t_1>t_0$,  the following Poincar\'e inequality holds:
$$   \EE \left[\left\{\EE\left(F(X_{[0,T]}^x)|\F_{t}\right)\right\}^2\right]-   \left\{\EE\left(F (X_{[0,T]}^x)\right)\right\}^2
 \le 2\int_0^t \E_{s,T}^{K,C}(F,F)\d s,~F\in \F C^{x,R,T}_{loc}.$$
\end{enumerate}\end{thm}

\paragraph{Remark 1.1.} $(a)$ We allow that the $L$-diffusion process is explosive in Theorem \ref{T1.1}.

$(b)$ When the function $K$ satisfy certain integrable conditions mentioned as in Theorem \ref{T1.2}, then Theorem \ref{T1.1} implies the associated global characterizations. In fact, we may derive easily the conclusion letting $R$ tends to $\infty$ in each term in Theorem \ref{T1.1}.

$(c)$ When $Z=0$, Theorem \ref{T1.1} first provides the characterizations for the upper bound of the Ricci curvature, these characterizations depend only on the local bounded property of Ricci curvature in Riemannian manifold. Our results are different in essence from characterizations of the global upper and lower bounds of the Ricci curvtaure. In fact, this is because that for each Riemannian manifold, the Ricci curvature is local bounded, but is not global bounded. A challenging problem: How to give a reasonable characterization of the upper bound of the Ricci curvature without any other geometrcial informations? After our work, Wang\cite{W17} and Sturm\cite{S} give some equivalent conditions about the upper bound of the Ricci curvature in Riemannian manifold and metric measure space respectively, but these equivalent conditions also depend on other geometrical informations.


The following Theorem \ref{T1.2} will characterize general lower and upper bounds of the Bakry-Emery curvature.

\beg{thm}\label{T1.2} Suppose $K_1$ and $K_2$ are two continuous functions on $M$ with $K_1\geq K_2$ such that
\beg{align} \label{DE} \EE\e^{(2+\vv)\int_0^T |K_1|(X_s^x)+|K_2|(X_s^x)\d s}<\infty \  {\rm for\ some}\ \vv,T>0.\end{align} For any $p, q\in [1,2]$, the following statements are equivalent each other:
\beg{enumerate}\item[$(1)$] For any $x\in M$,
\beg{align*} K_2(x)g_x\leq\Ric_Z(x)\le K_1(x)g_x.\end{align*}
\item[$(2)$]  For any $T>0$ and $x\in M$, $f\in C_0^\infty(M)$ with $|\nabla f|(x)=1$,
\beg{align*} &|\nn P_T f|^p(x) \le \EE\Big[\e^{-p\int^T_0 K_2(X_u^x)\d u }|\nn f|^p(X_T^x)\Big],\\
&\Big|\nn  f(x)- \ff 1 2\nn P_Tf(x) \Big|^q
 \le \EE\bigg[\Big(1+ \mu^{K_1,K_2}_{x,T}([0,T])\Big)^{q-1} \\
 &\times \bigg(\Big|\nn f(x)-\ff 1 2 A_T^{K_1,K_2}U_0^x(U_T^x)^{-1} \nn f(X_T^x)\Big|^q
+\ff { \mu^{K_1,K_2}_{x,T}([0,T])} {2^q}  \Big(A_T^{K_1,K_2}\Big)^q\big| \nn f(X_T^x)\big|^q  \bigg)\bigg].\end{align*}
\item[$(3)$] For any $F\in \F C_T^\infty, x\in M$ and $T>0$,
\beg{align*} |\nn_x \EE F(X_{[0,T]}^x)|^q &\le  \EE \bigg\{\Big(1+ \mu^{K_1,K_2}_{x,T}([0,T])\Big)^{q-1}\\
&\times\bigg( \big|\dot{D}^{K_1,K_2}_0F\big|^q +\int_0^T
  \Big|\dot{D}^{K_1,K_2}_sF(X_{[0,T]}^x)\Big|^q \mu^{K_1,K_2}_{x,T}(\d s)\bigg)\bigg\}.\end{align*}
\item[$(4)$] For any  $t_0,t_1\in [0,T]$ with $t_1>t_0$, and any $ x\in M$,  the following log-Sobolev inequality holds: $$\aligned
    &\EE\left[\EE\big(F^2(X_{[0,T]}^x)|\F_{t_1}\big) \log \EE(F^2(X_{[0,T]}^x)|\F_{t_1})\right]\\&-\EE\left[\EE
\big(F^2(X_{[0,T]}^x)|\F_{t_0}\big) \log \EE(F^2(X_{[0,T]}^x)|\F_{t_0})\right]
 \le 4 \int_{t_0}^{t_1} \E_{s,T}^{K_1,K_2}(F,F)\d s,~F\in \F C_T^\infty.\endaligned$$
 \item[$(5)$] For any  $t \in [0,T]$  and $ x\in M$,  the following Poincar\'e inequality holds:
$$   \EE \Big[\big\{\EE(F(X_{[0,T]}^x)|\F_{t})\big\}^2\Big]-   \Big\{\EE  \big[F (X_{[0,T]}^x)\big] \Big\}^2
 \le 2\int_0^t \E_{s,T}^{K_1,K_2}(F,F)\d s,\ \ F\in \F C_T^\infty.$$
\end{enumerate}\end{thm}

\begin{Remark} $(1)$ When $K_2=-K_1=-K$ for some constant $K\geq0$,  Naber\cite{N} first proved Theorem \ref{T1.2}.  After that, Naber's results had been extended to the case of manifold with a boundary by Wang-Wu\cite{WW} through another new method.

$(3)$ When we obtained Theorem \ref{T1.2}, we also noticed that Cheng-Thalmaier\cite{CT1} obtained some similar conclusions, their results are similar to (3), (4) and (5) of  our Theorem \ref{T1.2}, but the techniques are also somewhat different.
~Two papers were finished independently.

$(3)$ According to Theorem \ref{T1.2}, $M$ is an Einstein manifold with $\Ric=K$  if and only if all/some of items (2)-(5) hold for $K_1=K_2=K$ and $Z=0$.

\end{Remark}

The rest of this paper is organized as follows: In Section 2, we will present the proof of Theorem 1.2. As an application, some equivalent conditions for variable lower bounds of Ricci curvature will be given. Finally, the complete proof of Theorem 1.1 will be outlined in Section 3.

\section{Proof of Theorem \ref{T1.2}}
This section is devoted to mainly prove Theorem \ref{T1.2}.
To do that, we first introduce some basic notation partly coming from \cite{W2}. Let $f\in C_0^\infty(M)$ with $|\nn f(x)|=1$ and  $\Hess_f(x)=0$.  According to \cite[Theorem 3.2.3]{W2},  if $x\in M$ then for any $p>0$ we have
\beg{equation}\beg{split} \label{RIC}   \Ric_Z(\nn f,\nn f)(x)  &= \lim_{t\downarrow 0} \ff{P_t|\nn f|^p(x)-|\nn P_t f|^p(x)}{pt} \\
&= \lim_{t\downarrow 0} \ff 1 t \bigg(\ff{P_tf^2(x)-(P_tf)^2(x)}{2t}-|\nn P_tf(x)|^2\bigg).\end{split}\end{equation}
For any $s\ge 0$, consider the following resolvent equation
\begin{equation}\label{eq2.2}
\frac{\d Q^x_{s,t}}{\d t}=-Q^x_{s,t}\Ric_Z(U^x_t), \quad t\geq s,\ Q^x_{s,s}=\Id.\end{equation}
Then $(Q_{s,t}^x)_{t\ge s}$ is an adapted right-continuous process on $\R^d\otimes \R^d$.
For the sake of convenience, let $Q_t^x:=Q_{0,t}^x$. By \cite[Proposition 2.1]{Hsu0}(also refer to \cite[Lemma 4.2.3]{W2} and references within),
for any $F\in \F C_T^\infty$ with $F(\gg)= f(\gg_{t_1},\cdot, \gg_{t_N}), f\in C_0^{Lip}(M)$ and $0\le t_1<\cdots \leq t_N$,
\beq\label{eq2.3} (U_0^x)^{-1} \nn_x  \EE\big[F(X_{[0,T]}^x)\big] = \sum_{i=1}^N \EE\Big[Q_{t_i}^x(U_{t_i}^x)^{-1} \nn_i f(X_{t_1}^x,\cdots,X_{t_N}^x)\Big],\end{equation}
where $\nn_x$ denotes the gradient in $x\in M$ and $\nn_i$ is the gradient with respect to the $i$-th component.
In particular, taking $F(\gg)= f(\gg_t)$,  we obtain
\beq\label{GR2} \nn P_t f(x)= U_0^x \EE\big[Q_t^x(U_t^x)^{-1}\nn f(X_t^x)\big],\ \ x\in M, f\in C_0^\infty(M), t\ge 0.\end{equation}

Finally, for the above  $F\in \F C_T^\infty$, the damped gradient of $F$ is denoted by
\beq\label{TTD} \tt D_tF(X_{[0,T]}^x)= \sum_{i: t_i> t}Q_{t,t_i}^xU_{t_i}^{-1} \nn_if(X_{t_1}^x,\cdots, X_{t_N}^x),\ \
t\in [0,T].\end{equation}
Then  the martingale representation theorem in \cite{CHL}(see also \cite{W11, W2}) implies that
\beq\label{MF} \EE\big(F(X_{[0,T]}^x)\big|\F_t\big)= \EE[F(X_{[0,T]}^x)] +\ss 2 \int_0^t\big\<\EE(\tt D_sF(X_{[0,T]}^x) |\F_s), \d W_s\big\>,\ \ t\in [0,T].\end{equation}

\beg{proof}[Proof of Theorem \ref{T1.2}]
For simplicity, below we will 
write $F$ and $ f$ for  $F(X_{[0,T]}^x)$ and $f(X_{t_1}^x,\cdots, X_{t_N}^x)$ respectively. If $\Ric_Z$ has symmetric lower and upper bounds, A. Naber \cite{N} and Wang-Wu\cite{WW} obtained this conclusion for the constant bound and for pointwise symmetric bound respectively. But when the lower and upper bounds of $\Ric_Z$ is not symmetric, it is quite difficult to establish some functional inequalities by using the uniform norm of $\Ric_Z$ such that these inequalities can characterise the associated the lower and upper bounds. To overcome the difficulty, we may make a symmetrization of Bakry-Emery curvature such that it is symmetric, i.e. we may consider the curvature $$\Ric_Z^{K_1,K_2}(U^x_t):=\Ric_Z(U^x_t)-\frac{K_1(X_t^x)+K_2(X_t^x)}{2}\Id,$$
where $\Ric_Z(U^x_t):
\mathbb{R}^d\to\mathbb{R}^d$ is defined by
$$\<\Ric_Z(U^x_t)(a),b\>:=\Ric(U^x_ta,U^x_tb)-\<\nabla_{U^x_ta}Z,U^x_t b\>,\quad a,b\in\mathbb{R}^d.$$
Then by $(1)$, we get
\begin{equation}\label{eq2.7}\left\|\Ric_Z^{K_1,K_2}(U^x_t)\right\|\leq \frac{K_1(X_t^x)-K_2(X_t^x)}{2}.\end{equation}
Let $\widetilde{Q}^x_{s,t}$ be the solution of the following resolvent equation
\begin{equation}\label{eq2.8}
\frac{\d \widetilde{Q}^x_{s,t}}{\d t}=-\widetilde{Q}^x_{s,t}\Ric_Z^{K_1,K_2}(U^x_t),\quad t\geq s,\ \widetilde{Q}^x_{s,s}=\Id.
\end{equation}
and denote $\widetilde{Q}_t^x:=\widetilde{Q}^x_{0,t}$ for convenience.
Combining this with (\ref{eq2.8}),
\begin{equation}\label{eq2.9}
\widetilde{Q}^x_{s,t}=\e^{\frac{\int^t_s K_1(X_s^x)+K_2(X_s^x)\d s}{2}}Q^x_{s,t}.
\end{equation}
Following the line of \cite{WW}, in the following, we will present a full proof of Theorem \ref{T1.2}.

{\bf (a)}  (1) $\Rightarrow$ (3) for all $q\ge 1$. According to (\ref{eq2.3}), we have
\begin{equation*}U_0^{-1}\nabla_x\mathbb{E}[F]=\mathbb{E}\bigg[\sum_{i=1}^NQ^{x}_{t_i}(U_{t_i}^x)^{-1}\nabla_if\bigg].\end{equation*}
Then by (\ref{eq2.9}) and (\ref{eq1.7}),
\beg{equation}\label{eq2.10} \beg{split}
&U_0^{-1}\nabla_x\mathbb{E}[F]=\mathbb{E}\bigg[\sum_{i=1}^N\widetilde{Q}^{x}_{t_i}A_{t_i}^{K_1,K_2}(U_{t_i}^x)^{-1}\nabla_if\bigg]
 \\
&=\mathbb{E}\bigg[\sum_{i=1}^N\bigg(I+\int^{t_i}_0\widetilde{Q}^x_s\left(\Ric_Z(U^x_s)-\frac{K_1(X_s^x)+K_2(X_s^x)}{2}\Id\right)\d s
\bigg)A_{t_i}^{K_1,K_2}(U_{t_i}^x)^{-1}\nabla_if\bigg]\\
&=\mathbb{E}\bigg[\sum_{i=1}^NA_{t_i}^{K_1,K_2}(U_{t_i}^x)^{-1}\nabla_if\\
&~~~~~~~+\sum_{i=1}^N\bigg(\int^{t_i}_0\widetilde{Q}^x_s
\left(\Ric_Z(U^x_s)-\frac{K_1(X_s^x)+K_2(X_s^x)}{2}\Id\right)\d s\bigg)A_{t_i}^{K_1,K_2}(U_{t_i}^x)^{-1}\nabla_if\bigg]\\
\\&=\mathbb{E}\bigg[\dot{D}^{K_1,K_2}_0F+\int^T_0\bigg\{\widetilde{Q}^x_s\left(\Ric_Z(U^x_s)-\frac{K_1(X_s^x)+K_2(X_s^x)}{2}\Id\right)\dot{D}^{K_1,K_2}_sF\bigg\}\d s\bigg].\end{split} \end{equation}
In addition, by (\ref{eq2.7}) and (\ref{eq2.8}), we have
\beq\label{QQ2} \Big\|\widetilde{Q}^x_s\Big\|\leq \e^{\frac{\int^{s}_0 K_1(X_u^x)-K_2(X_u^x)\d u}{2}}.
\end{equation}  Combining these with (1), \eqref{MU} and using H\"older's inequality twice, we obtain
\beq\label{eq2.12}\aligned
&\big|\nabla_x\mathbb{E}[F]\big|^q \leq \bigg\{\EE \big|\dot{D}^{K_1,K_2}_0F\big|
+\EE \int^T_0 \big|\dot{D}^{K_1,K_2}_sF\big| \mu^{K_1,K_2}_{x,T}(\d s)\bigg\}^q\\
&\le  \EE\bigg\{ \big|\dot{D}^{K_1,K_2}_0F\big| + \int_0^T
  \big|\dot{D}^{K_1,K_2}_sF\big| \mu^{K_1,K_2}_{x,T}(\d s)\bigg\}^q\\
  &\le \EE \bigg\{\bigg(\big|\dot{D}^{K_1,K_2}_0F\big|^q +\ff {\big(\int_0^T
  \big|\dot{D}^{K_1,K_2}_sF\big| \mu^{K_1,K_2}_{x,T}(\d s) \big)^q} { \Big\{\mu^{K_1,K_2}_{x,T}([0,T]) \Big\}^{q-1}}  \bigg)\big(1+ \mu^{K_1,K_2}_{x,T}([0,T])\big)^{q-1}\bigg\}\\
  &\le \EE \bigg\{\bigg( \big|\dot{D}^{K_1,K_2}_0F\big|^q +\int_0^T
  \big|\dot{D}^{K_1,K_2}_sF\big|^q \mu^{K_1,K_2}_{x,T}(\d s)\bigg)\big(1+ \mu^{K_1,K_2}_{x,T}([0,T])\big)^{q-1}\bigg\}.\endaligned\end{equation}
Thus, (3) holds.

{\bf (b)} (3) $\Rightarrow$ (2) for all $p=q$. Taking $F(\gg)= f(\gg_T)$, then we have $\EE F(X_{[0,T]}^x)=P_T f(x)$ and
\beq\label{eq2.13}\dot{D}^{K_1,K_2}_s F=A_T^{K_1,K_2}(U_T^x)^{-1}\nabla f(X_T^x),\quad s<T.\end{equation}
Thus
 $\big|\dot{D}^{K_1,K_2}_s F\big|=A_T^{K_1,K_2}|\nn f(X_T)|$ for $s\in [0,T].$ Finally, (3) implies
the first inequality  in (2) with $p=q$. Similarly, by taking $F(\gg)=f(\gg_0)-\ff 1 2 f(\gg_T)$, we have
 $ \EE F= f(x) -\ff 1 2 P_T f(x)$ and
 \beq\label{eq2.14}\aligned
 & \left|\dot{D}^{K_1,K_2}_0 F\right|= \left|\nn f(x)-\ff 1 2 A_T^{K_1,K_2}U_0^x (U_T^x)^{-1}\nn f(X_T^x)\right|,\\
  &\left|\dot{D}^{K_1,K_2}_s F\right|\le \ff 1 2 A_T^{K_1,K_2}\left|\nn f(X_T^x)\right|,\ \ s\in (0,T].\endaligned\end{equation}
 Then the second inequality in (2) is also implied by (3).

{\bf (c)} (2) for some $p\ge 1$ and $q\in [1,2] \Rightarrow$ (1). Let $x\in M$. According to the first inequality in (2) and \eqref{RIC} yield
\beq\label{eq2.15} \beg{split} 0&\leq \lim_{T\rightarrow0}\frac{\EE\Big[\e^{-p\int^T_0 K_2(X_s^x)\d s }|\nn f|^p(X_T^x)\Big]- |\nn P_T f|^p(x)}{pT}\\
&=\lim_{T\rightarrow0}\EE\left[\frac{\e^{-p\int^T_0 K_2(X_s^x)\d s }-1}{pT}|\nn f|^p(X_T^x)\right] +\lim_{T\rightarrow0}\frac{P_T|\nn f|^p(x)- |\nn P_T f|^p(x)}{pT}\\
&=-K_2+\Ric_Z(\nn f,\nn f)\end{split}\end{equation}
 where $f\in C_0^\infty(M)$ with $\Hess_f(x)=0$ and $|\nn f(x)|=1$.
This implies  $\Ric_Z\ge K_2$.

Next, we prove that the second inequality in (2) implies $\Ric_Z\le K_1$. By H\"older's inequality,   the second inequality in (2) for some $q\in [1,2]$ implies the same
inequality  for $q=2$:
\beq\label{eq2.16} \beg{split} &\Big|\nn  f(x)- \ff 1 2\nn P_Tf(x) \Big|^2
 \le \EE\bigg[\big(1+ \mu^{K_1,K_2}_{x,T}([0,T])\big) \\
 &\times \bigg(\Big|\nn f(x)-\ff 1 2 A_T^{K_1,K_2}U_0^x(U_T^x)^{-1} \nn f(X_T^x)\Big|^2+\ff { \mu^{K_1,K_2}_{x,T}([0,T])} {4} (A_T^{K_1,K_2})^2\big| \nn f(X_T^x)\big|^2  \bigg)\bigg].\end{split}\end{equation}
Equivalently,
\beq\label{eq2.17} \beg{split}  &|\nn f(x)|^2-\<\nn f(x),\nn P_Tf(x)\>+\ff 1 4|\nn P_Tf(x)|^2\\
&\le \EE\bigg[\left(1+ \mu^{K_1,K_2}_{x,T}([0,T])\right) \times \bigg(|\nn f|^2(x)-\left\langle \nn f(x), A_T^{K_1,K_2}U_0^x(U_T^x)^{-1} \nn f(X_T^x)\right\rangle \\
 &~~~~~~~~~~~~~+\ff { \mu^{K_1,K_2}_{x,T}([0,T])+1} {4} \left(A_T^{K_1,K_2}\right)^2\big| \nn f\big|^2(X_T^x) \bigg)\bigg].\end{split}\end{equation}
Thus,
\beq\label{eq2.18} \beg{split}  &0\le \EE\bigg[\mu^{K_1,K_2}_{x,T}([0,T]) |\nn f|^2(x)\\&~~~~~~~+\langle \nn f(x), \nn P_Tf(x)-\big(1+ \mu^{K_1,K_2}_{x,T}([0,T])\big)A_T^{K_1,K_2}U_0^x(U_T^x)^{-1} \nn f(X_T^x)\rangle \\
 &~~~~~~~~+\ff 1 4 \left[\big(\mu^{K_1,K_2}_{x,T}([0,T])+1\big)^2(A_T^{K_1,K_2})^2\big| \nn f\big|^2(X_T^x)- |\nn P_Tf(x)|^2\right]\\&=
\EE\Big[\mu^{K_1,K_2}_{x,T}([0,T])\Big] |\nn f|^2(x)+\big\<\nn f(x), \nn P_Tf(x)- \EE \big[U_0^x (U_T^x)^{-1} \nn f(X_T^x)\big]\big\>
\\&~~~~~~~+\left\langle \nn f(x), \EE\Big[\big\lbrace1-\big(1+ \mu^{K_1,K_2}_{x,T}([0,T])\big)A_T^{K_1,K_2}\big\rbrace U_0^x(U_T^x)^{-1} \nn f(X_T^x)\Big]\right\rangle \\
 &~~~~~~~~+\ff 1 4 \EE\left[\Big\lbrace(\mu^{K_1,K_2}_{x,T}([0,T])+1)^2(A_T^{K_1,K_2})^2-1\Big\rbrace\big| \nn f\big|^2(X_T^x)\right]\\&
~~~~~~~~+\ff 1 4 \big[P_T|\nn f|^2(x)- |\nn P_Tf(x)|^2\big].\end{split}\end{equation}
By \eqref{eq2.3} and \eqref{GR2}, we have
\beq\label{eq2.19} \beg{split} &\big\<\nn f(x), \nn P_Tf(x)- \EE [U_0^x (U_T^x)^{-1} \nn f(X_T^x)]\big\>\\
&= -\int_0^T\big\<\nn f(x), U_0^x  \Ric_Z(U_r^x)(U_T^x)^{-1}\nn f(X_T^x)\big\>\d r\\
& =- T\Ric_Z(\nn f,\nn f)(x)+ {\rm o} (T).\end{split}\end{equation}
Combining the above these with (\ref{RIC}), we obtain
\beq\label{eq2.20} \beg{split}  &0\le \lim_{T\rightarrow0}\frac{\EE\Big[\mu^{K_1,K_2}_{x,T}([0,T])\Big]}{T}|\nn f(x)|^2\\
&~~~~+\lim_{T\rightarrow0}\frac{\big\<\nn f(x), \nn P_Tf(x)- \EE \big[U_0^x (U_T^x)^{-1} \nn f(X_T^x)\big]\big\>}{T}\\
&+\lim_{T\rightarrow0}\frac{\left\langle \nn f(x), \EE\Big[\Big\lbrace1-\Big(1+ \mu^{K_1,K_2}_{x,T}([0,T])\Big)A_T^{K_1,K_2}\Big\rbrace U_0^x(U_T^x)^{-1} \nn f(X_T^x)\Big]\right\rangle}{T}\\
&~~~~+\ff 1 4\lim_{T\rightarrow0}\frac{\EE\left[\Big\lbrace\left(\mu^{K_1,K_2}_{x,T}([0,T])+1\right)^2(A_T^{K_1,K_2})^2-1\Big\rbrace\big| \nn f\big|^2(X_T^x)\right]}{T}
\\
&~~~~+\ff 1 4\lim_{T\rightarrow0}\frac{\big[P_T|\nn f|^2(x)- |\nn P_Tf(x)|^2\big]}{T}\\
&=\frac{K_1-K_2}{2}|\nn f(x)|^2-\Ric_Z(\nn f,\nn f)(x)+K_2|\nn f(x)|^2
\\&~~~~~~~~-\frac{K_2}{2}|\nn f(x)|^2+\ff 1 2\Ric_Z(\nn f,\nn f)(x).\end{split}\end{equation}
This implies
$\Ric_Z(\nn f,\nn f)(x)\le K_1.$

{\bf (d)}  (5) $\Rightarrow$ (1). Let $F(\gg)= f(\gg_T)$. Then (5) implies
\beq\label{eq2.21} P_Tf^2(x)-  (P_Tf(x))^2 \le 2 \int^T_0\EE\Big[\Big(1+ \mu^{K_1,K_2}_{x,T}([s,T])\Big)^2\Big(A_{s,T}^{K_1,K_2}\Big)^2|\nn f|^2(X_T^x)\Big]\d s .\end{equation}
 For $f$ in \eqref{RIC}, combining this
with \eqref{RIC} we obtain
\beg{align*} &\Ric_Z(\nn f, \nn f)(x) \\
&=\lim_{T\to 0} \ff 1 T \bigg(\ff{P_T f^2(x)- (P_Tf)^2(x)}{2T}-|\nn P_T f|^2\bigg)\\
&\le \lim_{T\to 0} \ff 1 T \bigg(\ff{1}{T}\int^T_0\EE\Big[\Big(1+ \mu^{K_1,K_2}_{x,T}([s,T])\Big)^2\Big(A_{s,T}^{K_1,K_2}\Big)^2|\nn f|^2(X_T^x)\Big]\d s-|\nn P_T f|^2\bigg)
\\&=2\lim_{T\to 0}\frac{P_T|\nn f|^2(x)-|\nn P_Tf|^2(x)}{2T}\\&~~~~+|\nn f|^2(x)\lim_{T\to 0}\ff {\ff{1}{T}\int^T_0\EE\Big(A_{s,T}^{K_1,K_2}\Big)^2 \Big(1+ \mu^{K_1,K_2}_{x,T}([s,T])\Big)^2\d s-1} {T} \\
&\leq2\Ric_Z(\nn f, \nn f)(x)+|\nn f|^2(x)\\
&\times\lim_{T\to 0}\Bigg\{\ff {\ff{1}{T}\int^T_0\Big[\EE\Big(A_{s,T}^{K_1,K_2}\Big)^2-1\Big]\d s} {T}
+\ff {\ff{2}{T}\int^T_0\EE\Big[\Big(A_{s,T}^{K_1,K_2}\Big)^2\mu^{K_1,K_2}_{x,T}([s,T])\Big]\d s} {T}\Bigg\} \\
&=2\Ric_Z(\nn f, \nn f)(x)+\left(-K_1-K_2+K_1-K_2\right)|\nn f(x)|^2\\
&= 2\Ric_Z(\nn f,\nn f)(x) -2K_2|\nn f|^2(x).\end{align*}
This implies $\Ric_Z(\nn f,\nn f)(x)\ge K_2|\nn f(x)|^2$.

Next, we will prove the upper bound estimates. We take $F(\gg)= f(\gg_\vv)-\ff 1 2f(\gg_T)$ for $\vv\in (0,T).$  By \eqref{MG},
$$\aligned|\dot{D}^{K_1,K_2}_{t,s} F|&= \Big|A_{t,\vv}^{K_1,K_2}\nn f(X_\vv)-\ff 1 2 A_{t,T}^{K_1,K_2}U_\vv^x(U_T^x)^{-1}\nn f(X_T^x)\Big|1_{[0,\vv)}(s)\\
&+ \ff 1 2A_{t,T}^{K_1,K_2}|\nn f(X_T^x)|1_{[\vv,T]}(s).\endaligned$$
Then (5) implies
\beq\label{eq2.22}\beg{split} I_\vv &:= \EE\Big[f(X_\vv^x)-\ff 1 2 \EE(f(X_T^x)|\F_\vv)\Big]^2-\Big(P_\vv f(x)-\ff 1 2P_Tf(x)\Big)^2\\
&\le 2  \int^\vv_0\EE\bigg\{(1+ \mu^{K_1,K_2}_{x,T}([0,T]))\Big(\Big|A_{t,\vv}^{K_1,K_2}\nn f(X_\vv)-A_{t,T}^{K_1,K_2}\ff 1 2 U_\vv^x(U_T^x)^{-1}\nn f(X_T^x)\Big|^2 \\
&\qquad\quad \  +|\nn f(X_T^x)|^2\ff { \int^T_\vv\big(A_{\vv,T}^{K_1,K_2}\big)^2\mu^{K_1,K_2}_{x,T}(\d t)} 4\bigg)\bigg\}\d t+c \vv^2=: J_\vv,\ \ \vv\in (0,T)\end{split}\end{equation} for some constant $c>0$. It is easy to show that
\beq\label{eq2.23}\beg{split}  \lim_{\vv\to 0} \ff{J_\vv}\vv = \EE\bigg\{&\Big(1+ \mu^{K_1,K_2}_{x,T}([0,T])\Big)\Big(\Big|\nn f(x) - \ff 1 2A_T^{K_1,K_2} U_0^x(U_T^x)^{-1}\nn   f(X_T^x)\Big|^2\\
&+\big(A_{T}^{K_1,K_2}\big)^2\ff {\mu^{K_1,K_2}_{x,T}([0,T])|\nn f(X_T^x)|^2} 4\bigg)\bigg\}.\end{split}\end{equation}

On the other hand, we have
\beq\label{eq2.24} \beg{split} \ff{I_\vv}\vv =& \ff{P_\vv f^2-(P_\vv f)^2}\vv +\ff 1{4\vv} \EE\Big[\big\{\EE\big(f(X_T^x)|\F_\vv\big)\big\}^2 -(P_Tf)^2(x)\Big]\\
& + \ff{\EE[f(X_T^x)\{P_\vv f(x)-f(X_\vv^x)\}]}\vv.\end{split} \end{equation} Let $f\in C_0^\infty(M)$, we have
\beq\label{eq2.25} \lim_{\vv\to 0} \ff{P_\vv f^2-(P_\vv f)^2}\vv= 2|\nn f|^2(x). \end{equation}
Next, \eqref{TTD} and  \eqref{MF} yield
\beq\label{eq2.26} \EE(f(X_T^x)|\F_\vv)= P_Tf(x)+\ss 2 \int_0^\vv \big\<\EE\big(Q_{s,T}^x(U_T^x)^{-1}\nn f(X_T^x)\big|\F_s\big), \d W_s\big\>.\end{equation}
Then $$\EE[\EE(f(X_T^x)|\F_\vv)]^2= (P_Tf)^2 + 2 \int_0^\vv\EE |Q_{0,T}^x(U_T^x)^{-1}\nn f(X_T^x)|^2\d s.$$ This together with \eqref{GR2} leads to
\beq\label{eq2.27}\beg{split}  & \lim_{\vv\to 0}\ff 1{4\vv} \EE\Big[\big\{\EE\big(f(X_T^x)|\F_\vv\big)\big\}^2 -(P_Tf)^2(x)\Big]\\
 &=  \ff 1 2 \Big|\EE\big[Q_{0,T}^x(U_T^x)^{-1}\nn f(X_T^x)\big]\Big|^2 = \ff 1 2 |\nn P_T f(x)|^2.   \end{split}\end{equation}
 Finally, by It\^o's formula we have
\beg{align*} P_\vv f(x)-f(X_\vv^x)& = P_\vv f(x)-f(x) -\int_0^\vv Lf(X_s^x)\d s -\ss 2 \int_0^\vv \<\nn f(X_s^x), U_s^x\d W_s\>\\
&= {\rm o}(\vv)  -\ss 2 \int_0^\vv \<\nn f(X_s^x), U_s^x\d W_s\>.\end{align*} Combining this with \eqref{eq2.27} and \eqref{eq2.28}, we arrive at
$$\lim_{\vv\to 0} \ff{\EE[f(X_T^x)\{P_\vv f(x)-f(X_\vv^x)\}]}\vv=- 2\<\nn f(x),\nn P_t f(x)\>.$$
Substituting this and \eqref{eq2.27}-\eqref{eq2.29} into \eqref{eq2.26}, we obtain
$$\lim_{\vv\to 0}\ff{I_\vv}\vv= 2\Big|\nn f(x)- \ff 1 2 \nn P_Tf(x)\Big|^2.$$
Combining this with \eqref{eq2.24} and \eqref{eq2.25}, we prove the second inequality in (2) for $q=2$, which implies $\Ric_Z\le K_1$.

{\bf (f)} $(1)\Rightarrow (4).$ According to \eqref{MF},
\begin{equation}\label{eq2.28}G_t:= \EE(F^2|\F_t) =\mathbb{E}(F^2)+\sqrt{2}\int^t_0\big\<\EE(\tt D_s F^2|\F_s),\d W_s\big\>,\ \ t\in [0,T].\end{equation}
By It\^o's formula,
\begin{equation}\label{eq2.29}\aligned \d (G_t\log G_t)&=(1+\log G_t)\d G_t+\frac{|\EE(\tt D_s F^2|\F_s)|^2}{G_t}\d t\\
&\le  (1+\log G_t) \d G_t +4 \EE(|\tt D_sF|^2|\F_s)\d t.\endaligned\end{equation}
Then
\beq\label{eq2.30} \EE[G_{t_1}\log G_{t_1}]-\EE[G_{t_0}\log G_{t_0}]
\le 4 \int_{t_0}^{t_1}  \EE  |\tt D_s F|^2\d s .\end{equation}
By (\ref{TTD})  we have
$$\aligned&\tt D_sF=\sum_{i=1}^N1_{\{s< t_i\}}Q^x_{s,t_i}(U_{t_i}^x)^{-1}\nabla_if=\sum_{i=1}^N1_{\{s< t_i\}}A_{s,t_i}^{K_1,K_2}\widetilde{Q}^x_{s,t_i}(U_{t_i}^x)^{-1}\nabla_if
\\&=\sum_{i=1}^N1_{\{s< t_i\}}\bigg(I+\int^{t_i}_sQ^x_{s,t}\left\{ \Ric_Z(U_t^x)-\frac{K_1(X_t^x)+K_2(X_t^x)}{2}\right\} \bigg)A_{s,t_i}^{K_1,K_2}(U_{t_i}^x)^{-1}\nabla_if\d t
\\&=\dot{D}^{K_1,K_2}_{s,s}F+\int^T_sQ^x_{s,t} \left\{ \Ric_Z(U_t^x)-\frac{K_1(X_t^x)+K_2(X_t^x)}{2}\right\}\dot{D}^{K_1,K_2}_{s,t}F\d t.\endaligned$$
Combining this   with (1), \eqref{QQ2}, and using the Schwarz inequality, we prove
\begin{equation}\label{eq2.31}|\tt D_sF|^2 \le (1+\mu^{K,C}_{x,T}([s,T])) \bigg(|\dot{D}^{K_1,K_2}_{s,s}F|^2 +\int_s^T |\dot{D}^{K_1,K_2}_{s,t}F|^2\mu^{K,C}_{x,T}(\d t)\bigg).\end{equation}
This together with \eqref{eq2.30} implies the log-Sobolev inequality in (4).
\end{proof}

In particular, if $K_2,K_1\in\R$ are two constants with $K_1\geq K_2$, from Theorem \ref{T1.2} we easily obtain the following Corollary \ref{cor2.1}.

\beg{cor}\label{cor2.1} Suppose $K_2,K_1\in\R$ with $K_1\geq K_2$. For any $p, q\in [1,2]$, the following statements are equivalent each other:
\beg{enumerate}\item[$(1)$] For any $x\in M$,
\beg{align*} K_2g_x\leq\Ric_Z(x)\le K_1g_x.\end{align*}
\item[$(2)$]  For any $T>0$ and $x\in M$, $f\in C_0^\infty(M)$ with $|\nabla f|(x)=1$,
\beg{align*} &|\nn P_T f|^p(x) \le \e^{-K_2pT}P_T|\nn f|^p(x),\\
&\Big|\nn  f(x)- \ff 1 2   \nn P_Tf(x) \Big|^q
 \le \e^{\frac{K_1-K_2}{2}(q-1)T}\\
 & \times \EE \bigg[\Big|\nn f(x)-\ff 1 2 \e^{-\frac{K_1+K_2}{2}T}U_0^x(U_T^x)^{-1} \nn f(X_T^x)\Big|^q
  + \ff {\e^{\frac{K_1-K_2}{2}T}-1} {2^q}  \e^{-\frac{K_1+K_2}{2}qT}\big| \nn f(X_T^x)\big|^q\bigg].\end{align*}
\item[$(3)$] For any $F\in \F C_T^\infty, x\in M$ and $T>0$,
\beg{align*} |\nn_x \EE F(X_{[0,T]}^x)|^q &\le  \e^{\frac{K_1-K_2}{2}(q-1)T}\EE\bigg[|\dot{D}^{K_1,K_2}_0 F(X_{[0,T]}^x)|^q\\&
  ~~~~~~+ \int_0^T  |\dot{D}^{K_1,K_2}_s F(X_{[0,T]}^x)|^q\mu^{K,C}_{x,T}(\d s)\bigg].\end{align*}
\item[$(4)$] For any  $t_0,t_1\in [0,T]$ with $t_1>t_0$, and any $ x\in M$,  the following log-Sobolev inequality holds: $$\aligned
    &\EE\left[\EE\big(F^2(X_{[0,T]}^x)|\F_{t_1}\big) \log \EE(F^2(X_{[0,T]}^x)|\F_{t_1})\right]\\&-\EE\left[\EE
\big(F^2(X_{[0,T]}^x)|\F_{t_0}\big) \log \EE(F^2(X_{[0,T]}^x)|\F_{t_0})\right]
 \le 4 \int_{t_0}^{t_1} \E_{s,T}^{K_1,K_2}(F,F)\d s,\  \  F\in \F C_T^\infty.\endaligned$$
 \item[$(5)$] For any  $t \in [0,T]$  and $ x\in M$,  the following Poincar\'e inequality holds:
$$   \EE \Big[\big\{\EE(F(X_{[0,T]}^x)|\F_{t})\big\}^2\Big]-   \Big\{\EE  \big[F (X_{[0,T]}^x)\big] \Big\}^2
 \le 2\int_0^t \E_{s,T}^{K_1,K_2}(F,F)\d s,\ \ F\in \F C_T^\infty.$$
\end{enumerate}\end{cor}

\beg{cor}\label{cor2.2} Suppose $K$ is a continuous function on $M$ and
\beg{align} \label{eq2.32} \EE\e^{-p\int_0^T K^-(X_s^x)\d s}<\infty \  {\rm for\ some}\ \vv,T>0\end{align} for some $p\geq1$, the following two statements are equivalent:
\beg{enumerate}\item[$(1)$] For any $x\in M$,
\beg{align*}\Ric_Z(x)\geq K(x).\end{align*}
\item[$(2)$]  For any $T>0$ and $x\in M$, $f\in C_0^\infty(M)$,
\beg{align*} &|\nn P_T f|^p(x) \le \EE\Big[\e^{-p\int^T_0 K(X_u^x)\d u}|\nn f|^p(X_T^x)\Big].\end{align*}
\end{enumerate}\end{cor}

\beg{proof} Obviously, $(1)$ follows from $(2)$ by using (\ref{RIC}). In the following, it suffices to only show that $(1)\Rightarrow (2)$.
Let
\beg{align} \label{eq2.33}\aligned
&K_1(x)=\sup\left\{\Ric_Z(X,X):~X\in T_xM, ~|X|=1\right\}\\
&K_2(x)=\inf\left\{\Ric_Z(X,X):~X\in T_xM, ~|X|=1\right\}.
\endaligned\end{align}

(a) If $K_1$ and $K_2$ satisfy with the integrable condition (\ref{DE}), by Theorem \ref{T1.2}, then we have (3) of Theorem \ref{T1.2}. In particular, for the function $F(\gamma):=f(\gamma_T)$, it implies that
\beg{align} \label{eq2.34}\aligned
&|\nn P_T f|(x) \le \EE\Big[\e^{-\int^T_0 K_2(X_u^x)\d u}|\nn f|(X_T^x)\Big]
\endaligned\end{align}
Thus, according to (\ref{eq2.32}) and Cauchy-Schwartz inequality, (2) are implied by $K\leq K_2$.

(b) In general, by the similar argument in Section 3, we construct a sequence of processes $X^{x,R}_\cdot$ such that $X^{x,R}_t=X^x_t$ $\P_x$-$a.s.$ for every $t \le \tau_R$ on $M_R$.
Let $\Ric^{(R)}, \nabla^{(R)}$ and $\|\cdot\|_\varphi$ be the associated Ricci
curvature, the Levi-Civita connection, and the norm of vectors on
$M_\varphi$ respectively. Define $K_1^R, K_2^R$ as (\ref{eq2.32}) with $\Ric_Z$
replaced by $\Ric^{(R)}_{Z^R}$. Then $K_1^R, K_2^R$ satisfy with (\ref{DE}), by (a),
\beg{align} \label{eq2.35}\aligned
&|\nn \EE[f(X^{x,R}_T)] |^p \le \EE\Big[\e^{-p\int^T_0 K^R_2(X_u^x)\d u}|\nn f|^p(X_T^x)\Big]
\endaligned\end{align}
By (\ref{eq2.32}), (1) and the dominated convergence theorem, we have
\beg{align} \label{eq2.36}\aligned
\lim_{R\rightarrow\infty}\EE\Big[\e^{-p\int^T_0 K^R_2(X_u^x)\d u}|\nn f|^p(X_T^x)\Big]&=\EE\Big[\e^{-p\int^T_0 K_2(X_u^x)\d u}|\nn f|^p(X_T^x)\Big]\\&\leq \EE\Big[\e^{-p\int^T_0 K(X_u^x)\d u}|\nn f|^p(X_T^x)\Big]
\endaligned\end{align}
and
\beg{align} \label{eq2.37}\aligned
\lim_{R\rightarrow\infty}|\nn \EE[f(X^{x,R}_T)] |^p&=\lim_{R\rightarrow\infty}|\EE[Q_T^{R,*}U_T^R\nabla f(X^{x,R}_T)] |^p
\\&=|\EE[Q_T^*U_T^R\nabla f(X^x_T)] |^p=|\nn \EE[f(X^x_T)] |^p.\endaligned\end{align}
Applying (\ref{eq2.36}) and (\ref{eq2.37}) into (\ref{eq2.35}), this concludes the proof of Corollary \ref{cor2.2}.
\end{proof}

\begin{Remark}
The equivalence of Corollary \ref{cor2.2} was first proved by \cite{W2} under a different condition.
\end{Remark}

\section{Proof of Theorem \ref{T1.1}}
In this section, our aim is to prove Theorem \ref{T1.1}. To do that, we need to make some preparations. Fix $x\in M$ and for any $R>0$, the stopping time is given by
$$\tau_R=\inf\left\{t\ge 0: \rr(x,X_t^x)\ge R\right\}.$$ By \cite[Lemma 3.1.1]{W2} (see also \cite[Lemma 2.3]{ATW09}), there exists a constant $c>0$ such that
\beq\label{LO} \P(\tau_R\le T) \le \e^{-c/T},\ \ T\in (0,1].\end{equation}
For each $T>0$, let $\{t_i\}^\infty_{i=1}$ be the countable dense subset of the interval $[0,T]$. For any $l\in C_0^\infty(\R)$, we construct a family of functions as follows:
\begin{equation}\label{c*}
\rho_x^m(\gamma):=\sup_{1\le i \le m}d_M(\gamma(t_i),x),\ \
 m\geq1, \gamma \in W_x(M)
\end{equation}
and
$$l_x^m(\gamma):=l(\rho_x^m(\gamma)),\quad m\geq1,\gamma \in W_x(M).$$
It is obvious that $l_x^m:=l(\rho_x^m) \in \F C^\infty_b$.
Define
$$g(s):=\max_{1\le i \le m}s_i, \quad s=(s_1,\dots,s_m)\in \R^m.$$
Then
$$\rho_x^m=g\Big(\big(d_M(\gamma(t_1),x),\dots,d_M(\gamma(t_m),x)\big)\Big).$$
Since $g(s)$ is a Lipschitz continuous function on $\R^m$
with Lipschitz constant $1$, and $|\nabla d_M(\cdot,x)|\le 1$, then for every $m \ge 1$, we have
\begin{equation*}
\|D^{K_1,K_2} l_x^m(\gamma)\|_{\H}\le \sup_{r \in \R}|l'(r)|\e^{\int_0^T\frac{|K_1(\gamma_s)+K_2(\gamma_s)|}{2}\d s},\ \ \ \P_x-a.s. \ \gamma \in W_x(M).
\end{equation*}
Hence by Mazur's theorem (or refer to the argument in the proof of \cite[Lemma 2.2]{A}, \cite[Proposition 3.1]{RS} or \cite[Lemma 2.1]{CW}), we know
\begin{equation}\label{eq3.3}
 \|D^{K_1,K_2}l(\rho_x(\gamma))\|_{\H}\le \sup_{r \in \R}|l'(r)|\e^{\int_0^T\frac{|K_1(\gamma_s)+K_2(\gamma_s)|}{2}\d s},\ \ \ \P_x-a.s. \ \gamma \in W_x(M).
 \end{equation}
In addition, by the standard procedure and the local integration by parts formula(refer to the proof of Theorem 1.1 in \cite{CW}), we have 

 \begin{equation}\label{eq3.4}
 l(\rho_x)=\lim_{m\rightarrow\infty}l_x^m, \quad \E_1-\text{norm}.
 \end{equation}

\begin{lem}\label{l3.1}
Assume that $l \in C_0^{\infty}(\R)$ with $\supp (l)\subset B_{2R}(x),~l|_{B_{R}(x)}=1$ for some $R>0$. Then for any $p>0$ and for any $f\in C_0^\infty(M)$,
 \begin{equation}\label{eq3.5}
 \lim_{T\rightarrow\infty}\frac{1}{T^p}\nn_x \EE \Big[\big(1-l\left(\rho_x\right)\big)f(X^x_T)\Big]=0.
 \end{equation}
\end{lem}

\beg{proof} It suffices to only show that
\begin{equation}\label{eq3.6}
 \lim_{T\rightarrow\infty}\frac{1}{T^p}\d_x\Big(\EE \Big[\big(1-l\left(\rho_x\right)\big)f(X^x_T)\Big]\Big)(V)=0,\quad \forall~V\in T_xM.
 \end{equation}
By Lemma 2.1 in \cite{CLW1}, we know that there exists a $L^2$-integrable cut-off function $h:[0,T]\times W_x(M)\rightarrow \R^n$ such that
$$
h(\cdot,\gamma)=0,\quad \forall~\gamma\cap B_R(x)\neq\emptyset.$$
Then by the same argument of Lemma 3.8 in \cite{CLW2},
$$\EE \Big[\big(1-l\left(\rho_x\right)\big)f(X^x_T)\Big]$$
is differentiable, i.e.
\begin{equation}
\label{df}\beg{split}
&\d_x\Big(\EE \Big[\big(1-l\left(\rho_x\right)\big)f(X^x_T)\Big]\Big)(V(x))=\EE\Big(\d\Big[\big(1-l\left(\rho_x\right)\big)f(X^x_T)\Big](U_\cdot(x) h^V(x)\Big)
\\
&-\EE_x  \Big(\Big[\big(1-l\left(\rho_x\right)\big)f(X^x_T)\Big]\int_0^T\left \<\dot h^V_r+\frac{1}{2}U_r^{-1}(\sigma)\Ric_Z^{\#}(h^V_r), \d B_r\right\> \Big).
\end{split} \end{equation}
where $h^V(x)= h_\cdot(x)-h_0(x)+U_0^{-1}V$ and $B_s$ is  the stochastic anti-development of the $L$-diffusion processes. Thus,
by (\ref{LO}), we have
\begin{equation}\label{eq3.8}\beg{split}
 &\Big|\d_x\Big(\EE \Big[\big(1-l\left(\rho_x\right)\big)f(X^x_T)\Big]\Big)(V(x))\Big|\leq \EE\Big|\d\Big[\big(1-l\left(\rho_x\right)\big)f(X^x_T)\Big](U_\cdot(x) h^V(x)\Big|
\\
&+\EE_x  \bigg|\Big[\big(1-l\left(\rho_x\right)\big)f(X^x_T)\Big]\int_0^T\left \<\dot h^V_r+\frac{1}{2}U_r^{-1}(\sigma)\Ric_Z^{\#}(h^V_r), \d B_r\right\> \bigg|\\& \leq C_1\P(\tau_R\le T) \le C_1\e^{-C_2/T}
\end{split} \end{equation}
for some constants $C_1,C_2>0$.
This implies (\ref{eq3.6}).
\end{proof}

To prove Theorem \ref{T1.1}, the integration by parts formula with respect to Wiener measure on $W(M)$ is crucial. Since Theorem \ref{T1.1} does not require that the diffusion process generated by $L$ is non-explosive, that is to say, there are no any integrable conditions about the Ricci curvature. This means that the overall integration by parts formula does not hold. To overcome the difficulty, we may obtain the local the formula of integration by parts
by using the cutoff method. The idea is to  make a conformal change of metric such that the new Riemannian manifold is with bounded curvature(see \cite{WW,CW,WW2}) and two metrics are the same in a compact set. In fact, for any $R>0$,
taking
$\varphi\in C^\infty_0(M)$ such that $\varphi|_{B_{R+1}(x)}=1$. Let
$$\tau_{R}:=\inf\Big\{t\geq0: \rho_x(X^x_t)\geq R\Big\}$$
and
$$M_R:=\{y\in M: \varphi(y)>0\}.$$
According to \cite[section 2]{TW}(See also \cite{TW,FWW,WW2, CW}) and references in, $(M_R, \<\cdot, \cdot\>_R)$ is a
complete Riemannian manifold under the metric $$\<\cdot, \cdot\>_R:=\varphi^{-2}\<\cdot, \cdot\>, $$and
$$L_R=\varphi^2L=\frac{1}{2}\triangle^{(R)}+Z^{(R)}$$
for $\triangle^{(R)}$ the Laplace operator on $M_R$ and $Z^{(R)}$
some vector field on $M_\varphi$ such that
$$\displaystyle\sup_{M_R}\big(\|\Ric^{(R)}\|_\varphi+\|\nabla^{(R)}Z^{(n)}\|_\varphi\big)<\infty,$$
where $\Ric^{(R)}, \nabla^{(R)}$ and $\|\cdot\|_R$ are the Ricci
curvature, the Levi-Civita connection, and the norm of vectors on
$M_R$ respectively. Therefore, letting $\mathbb{P}^{x,R}$ be the
distribution of the $L_R$-diffusion process $X^{x,R}_\cdot$ on $M_R$. Then, we have the Driver's formula (\ref{eq3.9}) and the martingale representation theorem (\ref{MF}), but where the process $X^x_\cdot$ is now replaced by $X^{x,R}_\cdot$. Then, by repeating the precious computations of \eqref{eq2.10}, we consider the gradient formula on the new path space $W_x^T(M)$: for any $F\in \F C_T^\infty$, we have
\beg{equation*}\beg{split}
&U_0^{-1,R}\nabla_x\mathbb{E}_{\P^{x,R}}[F]=\mathbb{E}\bigg[D^{K_1^R,K_2^R}_0F
\\&~~~~~~~~~~~~~~~~~~~~~+\int^T_0\bigg\{\widetilde{Q}^{x,R}_s\bigg(\Ric^R_Z(U^{x,R}_s)-\frac{K^R_1(X_t^{x,R})+K^R_2(X_t^{x,R})}{2}\Id\bigg)\dot{D}^{K_1^R,K^R_2}_sF\bigg\}\d s\bigg],\end{split}\end{equation*}
where $U_0^{-1,R},D^{K_1^R,K_2^R}_0,\Ric^R_Z$ are defined similar to the previous one. 
In addition, by the above construction, we know $X^{x,R}_t=X^x_t$ $\P_x$-$a.s.$ for every
$t \le \tau_R$. Then for any $\gamma\in W_x^T(M)$ with $\gamma\subset B_R(x)$, then $\gamma$ may be looked as a path of  $W_x^T(M_R)$ and 
\beg{equation}\label{eq3.9}U_0^{-1}(\gamma)=U_0^{-1,R}(\gamma),~D^{K_1,K_2}_0F(\gamma)=D^{K_1^R,K_2^R}_0F(\gamma),~\Ric^R_Z(U^{x}_\cdot)=\Ric^R_Z(U^{x,R}_\cdot)\end{equation}
for all $F\in \F C_T^\infty.$
Thus, by \eqref{eq3.4}, lemma \ref{l3.1} and the dominated convergence theorem, for any $F\in \F C^\infty_b$ we have
\beg{equation*}\beg{split}
&U_0^{-1}\nabla_x\mathbb{E}_{\P^x}[l(\rho_x)F]=U_0^{-1}\lim_{m\rightarrow\infty}\nabla_x\mathbb{E}_{\P^x}[l(\rho^m_x)F]\\&=\lim_{m\rightarrow\infty}\mathbb{E}\bigg[D^{K_1^R,K_2^R}_0[l(\rho^m_x)F]
\\&~~~~~~~~~~~~~~~~~~~+\int^T_0\bigg\{\widetilde{Q}^{x,R}_s\bigg(\Ric^R_Z(U^{x,R}_s)-\frac{K^R_1(X_t^{x,R})+K^R_2(X_t^{x,R})}{2}\Id\bigg)\dot{D}^{K_1^R,K^R_2}_s[l(\rho^m_x)F]\bigg\}\d s\bigg]
\\&=\mathbb{E}\bigg[D^{K_1^R,K_2^R}_0[l(\rho_x)F]
\\&~~~~~~~~~~~~~~~~~~~+\int^T_0\bigg\{\widetilde{Q}^{x,R}_s\bigg(\Ric^R_Z(U^{x,R}_s)-\frac{K^R_1(X_t^{x,R})+K^R_2(X_t^{x,R})}{2}\Id\bigg)\dot{D}^{K_1^R,K^R_2}_s[l(\rho_x)F]\bigg\}\d s\bigg]
.\end{split} \end{equation*}
In particular, taking suitable functions $K_1,K_2$ such that $K_2\big|_{B_R(x)}=C, K_1\big|_{B_R(x)}=K(x)$, then the above equality implies that
\beg{equation}\label{eq3.10} \beg{split}
&U_0^{-1}\nabla_x\mathbb{E}_{\P^x}[l(\rho_x)F]=\mathbb{E}\bigg[D^{K,C}_0[l(\rho_x)F]\\&+\int^T_0\bigg\{\widetilde{Q}^{x}_s\bigg(\Ric^R_Z(U^{x}_s)
-\frac{K_1(X_t^{x})+C}{2}\Id\bigg)D^{K,C}_s[l(\rho_x)F]\bigg\}\d s\bigg],\quad F\in \F C^\infty_b. \end{split} \end{equation}

By the similar argument, we obtain
\beq\label{LMF}\beg{split} &\EE\big(F(X_{[0,T]}^x)\big|\F_t\big)\\&= \EE[F(X_{[0,T]}^x)] +\ss 2 \int_0^t\big\<\EE(\tt D_sF(X_{[0,T]}^x) |\F_s), \d W_s\big\>,\ \ t\in [0,T], F\in \F C^{x,R,T}_{loc}.\end{split}\end{equation}

\beg{proof}[Proof of Theorem \ref{T1.1}] By repeating the previous part of the proof of Theorem \ref{T1.2} and using (\ref{eq3.10}) and (\ref{LMF}), then (2)-(5) are implied by (1). Conversely, it is obvious that (2)$\Rightarrow$ (3) and (4)$\Rightarrow$ (5). Thus it suffices to show that (3)$\Rightarrow$ (1) and (5)$\Rightarrow$ (1).

(a) (3)$\Rightarrow$ (1): Let $F\in \F C^{x,R,T}_{loc}$ with $F(\gamma):=l\left(\rho_x\left(\gamma\right)\right) f(x)-\frac{1}{2}l\left(\rho_x\left(\gamma\right)\right)f(\gamma_T),$ where $l\in C_0^\infty(\R),f\in C_0^\infty(M)$ and $\supp (l)\subset B_R(x)$ and $l|_{B_{R/2}(x)}=1$. And satisfying with the following conditions
\beq\label{eq3.12} \beg{split}
&\supp (l)\subset B_R(x),\quad l|_{B_{R/2}(x)}=1\\&
f(x)=0,\quad |\nn f(x)|=1,\quad \Hess_f(x)=0.
\end{split}\end{equation}
By (\ref{LO}), we have
\beq\label{eq3.13} \beg{split} &\Big|\EE\big[\big(l\left(\rho_x\right)-1\big)\nn  f(x)\big]\Big|^2\leq \|\nn  f\|^2_\infty\P(\tau_{R/2}\leq T)=o(T^3).
\end{split}\end{equation}
Combining the above inequality with Lemma \ref{l3.1}, we have
\beq\label{eq3.14} \beg{split} &\left|\EE\big[l\left(\rho_x\right)\nn  f(x)\big]- \ff 1 2\nn \EE \big[l\left(\rho_x\right)f(X^x_T)\big]\right|^2\\
&=\left|\nn  f(x)- \ff 1 2\nn P_Tf(x)+\EE\big[(l\left(\rho_x\right)-1)\nn  f(x)\big]- \ff 1 2\nn \EE \big[(l\left(\rho_x\right)-1)f(X^x_T)\big]\right|^2\\
&=\left|\nn  f(x)- \ff 1 2\nn P_Tf(x)\right|^2+o(T^3).
\end{split}\end{equation}
According to (3) and (\ref{eq3.13}),
\beq\label{eq3.15} \beg{split} &\Big|\EE\big[l\left(\rho_x\right)\nn  f(x)\big]- \ff 1 2\nn \EE \big[l\left(\rho_x\right)f(X^x_T)\big]\Big|^2
=\big|\nabla\EE\big[F(X_{[0,T]}^x)\big]\big|^2\\&=\Big|\nn  f(x)\EE\big[l\left(\rho_x\right)\big]- \ff 1 2\nn \EE \big[l\left(\rho_x\right)f(X^x_T)\big]\Big|^2
\\&
 \le \EE\bigg[\big(1+ \mu^{K,C}_{x,T}([0,T])\big)\times \bigg(\Big| D^{K,C}_0F(X_{[0,T]}^x)\Big|^2+\int_0^T
  \Big|\dot{D}^{K,C}_sF(X_{[0,T]}^x)\Big|^2 \mu^{K,C}_{x,T}(\d s)\bigg)\bigg]
 \\&
 \le \EE\bigg[\big(1+ \mu^{K,C}_{x,T}([0,T])\big)\times \bigg(\Big| D^{K,C}_0\big[f(X^x_0)l\left(\rho_x\right)\big]-\ff 1 2 D^{K,C}_0\big[f(X_T^x)l\left(\rho_x\right)\big]\Big|^2\\
 &~~~~~~~~~~~~~~~~~~~+\ff 1 4\int_0^T
 \Big|\dot{D}^{K,C}_s\big[f(X_T^x)l\left(\rho_x\right)\big]\Big|^2 \mu^{K,C}_{x,T}(\d s)\bigg)\bigg]\\&
 \le \EE\bigg[\big(1+ \mu^{K,C}_{x,T}([0,T])\big)\times \bigg(\Big|l\left(\rho_x\right)\Big[\nn f(x)-\ff 1 2 A_T^{K,C}U_0^x(U_T^x)^{-1} \nn f(X_T^x)\Big]
 \\&~~~~-\ff 1 2 f(X_T^x)l'\left(\rho_x\right)D^{K,C}_s\rho_x\Big|^2+\\&~~~~\ff 1 4\int_0^T
 \Big|l\left(\rho_x\right)A_T^{K_1,C_R}U_0^x(U_T^x)^{-1} \nn f(X_T^x)+f(X_T^x)l'\left(\rho_x\right)D^{K,C}_s\rho_x\Big|^2 \mu^{K,C}_{x,T}(\d s)\bigg)\bigg]
\\&=\EE\bigg[\big(1+ \mu^{K,C}_{x,T}([0,T])\big)\times \bigg(
\Big|l\left(\rho_x\right)\nn f(x)- \ff 1 2 A_T^{K,C}U_0^x(U_T^x)^{-1} l\left(\rho_x\right)\nn f(X_T^x)\Big|^2\\&
~~~~+\ff {\mu^{K,C}_{x,T}([0,T])}{4} \big|A_T^{K,C}\big|^2\big|l\left(\rho_x\right)\nn f(X_T^x)\big|^2\bigg)\bigg]+C\P(\tau_{R/2}\leq T)\\&=\EE\bigg[\big(1+ \mu^{K,C}_{x,T}([0,T])\big)\times \bigg(
\Big|l\left(\rho_x\right)\nn f(x)- \ff 1 2 A_T^{K,C}U_0^x(U_T^x)^{-1} l\left(\rho_x\right)\nn f(X_T^x)\Big|^2\\&
~~~~+\ff {\mu^{K,C}_{x,T}([0,T])}{4} \big|A_T^{K,C}\big|^2\big|l\left(\rho_x\right)\nn f(X_T^x)\big|^2\bigg)\bigg]+o(T^3)
\\&=\EE\bigg[\big(1+ \mu^{K,C}_{x,T}([0,T])\big)\times \bigg(
\Big|\nn f(x)- \ff 1 2 A_T^{K,C}U_0^x(U_T^x)^{-1} \nn f(X_T^x)\Big|^2\\&
~~~~+\ff {\mu^{K,C}_{x,T}([0,T])}{4} \big|A_T^{K,C}\big|^2\big|\nn f(X_T^x)\big|^2\bigg)\bigg]+o(T^3).\end{split}\end{equation}
Next, using the same argument of (2)$\Rightarrow$ (1) in the proof of Theorem \ref{T1.2}. we obtain (1).

(a) (5)$\Rightarrow$ (1). Taking $f$ and $l$ as in (\ref{eq3.12}), but $$F(\gamma):=l\left(\rho_x\left(\gamma\right)\right) f(\gamma_\varepsilon)-\frac{1}{2}l\left(\rho_x\left(\gamma\right)\right)f(\gamma_T), \quad\varepsilon>0,$$ then  $F\in \F C^{x,R,T}_{loc}$ and by (\ref{MG})
$$\aligned|\dot{D}^{K,C}_{t,s} F|&= \Big|A_{t,\vv}^{K,C}\nn f(X_\vv)-\ff 1 2 A_{t,T}^{K,C}U_\vv^x(U_T^x)^{-1}\nn f(X_T^x)\Big|1_{[0,\vv)}(s)\\
&+ \ff 1 2A_{t,T}^{K,C}|\nn f(X_T^x)|1_{[\vv,T]}(s).\endaligned$$
Moreover, by (5) and \eqref{LO}, there exists a constant $C_1>0$ depending on $f,l,C,R$ and $K$ such that for any $\vv,T\in (0,1)$,
\beq\label{eq3.16}\aligned
I_\vv&:= \EE\Big[ \EE\Big(l\left(\rho_x\left(X_{[0,T]}^x\right)\right)f(X_\vv^x)-\ff 1 2l\left(\rho_x\left(X_{[0,T]}^x\right)\right)f(X_T^x)\Big|\F_\vv\Big)\Big]^2\\
&~~~~~~~~~~~~~~~~~-\Big[ \EE\Big(l\left(\rho_x\left(X_{[0,T]}^x\right)\right)f(X_\vv^x)-\ff 1 2l\left(\rho_x\left(X_{[0,T]}^x\right)\right)f(X_T^x)\Big)\Big]^2\\&
\le 2\int^\varepsilon_0\EE\bigg\{\big(1+\mu^{K,C}_{x,T}([t,T])\big)\bigg(
   |l\left(\rho_x\left(X_{[0,T]}^x\right)\right)\dot D_{t,t}^{K,C} F|^2
\\&~~~~~~~~~~~~~~~~ +\int^{T}_{t}   |l\left(\rho_x\left(X_{[0,T]}^x\right)\right)\dot D_{t,s}^{K,C} F |^2\mu^{K,C}_{x,T}(\d s)\bigg)\bigg\}\d t+ C_1\vv T^4.
\endaligned\end{equation}
Then
\beq\label{eq3.17}\beg{split}  \limsup_{\vv\to 0} \ff{I_\vv}\vv \le  \EE\bigg\{&l\left(\rho_x\left(X_{[0,T]}^x\right)\right)(1+\mu^{K,C}_{x,T}([0,T]))\Big(\Big|\nn f(x) - \ff 1 2 A_T^{K,C}U_0^x(U_T^x)^{-1}\nn   f(X_T^x)\Big|^2\\
&+\ff {l\left(\rho_x\left(X_{[0,T]}^x\right)\right)\mu^{K,C}_{x,T}([0,T])}4 |A_T^{K,C}\nn f|^2(X_T^x)\Big)\bigg\}+ {\rm o}(T^3) \end{split}\end{equation} for small $T>0$.
According to (\ref{eq2.24}), we have
\beq\label{eq3.18} \beg{split} \ff{I_\vv}\vv =& \ff{P_\vv f^2-(P_\vv f)^2}\vv +\ff 1{4\vv} \EE\Big[\big\{\EE\big(f(X_T^x)|\F_\vv\big)\big\}^2 -(P_Tf)^2(x)\Big]\\
& + \ff{\EE[f(X_T^x)\{P_\vv f(x)-f(X_\vv^x)\}]}\vv+o(T^3)\\
&= 2\Big|\nn f(x)- \ff 1 2 \nn P_Tf(x)\Big|^2+o(T^3). \end{split}\end{equation}
Combining this with  (\ref{eq3.16}),  we get
\beq\label{eq3.19}\beg{split}
&2\Big|\nn f(x)- \ff 1 2 \nn P_Tf(x)\Big|^2\\
&\leq\EE\bigg\{l\left(\rho_x\left(X_{[0,T]}^x\right)\right)(1+\mu^{K,C}_{x,T}([0,T]))\Big(\Big|\nn f(x) - \ff 1 2 U_0^x(U_T^x)^{-1}A_T^{K,C}\nn   f(X_T^x)\Big|^2\\
&+\ff {l\left(\rho_x\left(X_{[0,T]}^x\right)\right)\mu^{K,C}_{x,T}([0,T])}4 |A_T^{K,C}\nn f|^2(X_T^x)\Big)\bigg\}+o(T^3)\\
&=\EE\bigg\{(1+\mu^{K,C}_{x,T}([0,T]))\Big(\Big|\nn f(x) - \ff 1 2 U_0^x(U_T^x)^{-1}A_T^{K,C}\nn   f(X_T^x)\Big|^2\\
&+\ff {\mu^{K,C}_{x,T}([0,T])}4 |A_T^{K,C}\nn f|^2(X_T^x)\Big)\bigg\}+o(T^3)\end{split}
 \end{equation}
Using this estimate, the rest of the argument is similar  to the last part in the proof of (2) $\Rightarrow$ (1) of Theorem \ref{T1.2}.
\end{proof}

\section*{Acknowledgments}
It is a pleasure to thank M. Ledoux, F.Y. Wang and X. Chen for useful conversations and also thank L.J. Cheng and A. Thalmaier presented the draft of their paper\cite{CT1}. This research is supported by NNSFC (11371099)

\beg{thebibliography}{99}

\leftskip=-2mm
\parskip=-1mm





\bibitem{ATW09}  M. Arnaudon, A. Thalmaier, F.-Y. Wang, \emph{Gradient estimates and Harnack inequalities on non-compact Riemannian manifolds,} Stoch. Proc. Appl. 119(2009), 3653--3670.


\bibitem{Be} A. L. Besse, \emph{Einstein Manifolds,} Springer,
Berlin, 1987.


\bibitem{CHL} B. Capitaine,  E. P. Hsu,  M. Ledoux, \emph{Martingale representation and a simple proof of logarithmic
Sobolev inequalities on path spaces,} Elect. Comm. Probab. 2(1997), 71--81.

\bibitem{CW} X. Chen, B. Wu, \emph{Functional inequality on path space over a non-compact Riemannian manifold,}  J.
Funct. Anal. 266(2014), 6753-6779.

\bibitem{CLW1} X. Chen, X. M. Li, B. Wu, \emph{Quasi-regular Dirichlet form on loop spaces,}  Preprint.

\bibitem{CLW2} X. Chen, X. M. Li, B. Wu, \emph{Analysis on Free Riemannian Loop Space,}  Preprint.

\bibitem{CT1} L.J. Cheng, A. Thalmaier, \emph{Characterization of pinched Ricci curvature by functional inequalities,} {\it arXiv: 1611.02160}.

\bibitem{CT2} J. Cheng, A. Thalmaier, \emph{Spectral gap on Riemannian path space over static and evolving manifolds,} J. Funct. Anal. 274(2018) 959¡V984
\bibitem{D} B. Driver,  \emph{A Cameron-Martin type quasi-invariant theorem for Brownian motion on a compact Riemannian manifold},  J. Funct.
Anal. 110(1992), 272--376.

\bibitem{EL} K. D. Elworthy, X.- M. Li and Y. Lejan, \emph{On The geometry of diffusion
operators and stochastic flows,} Lecture Notes in Mathematics, 1720(1999), Springer-Verlag.


\bibitem{ES} O. Enchev and D. W. Stroock, \emph{Towards a Riemannian geometry on the path space over a Riemannian manifold,}
 J. Funct. Anal. 134 : 2 (1995), 392--416.

\bibitem{F} S. Fang, \emph{In\'egalit\'e du type de Poincar\'e sur l'espace
des chemins riemanniens}, C.R. Acad. Sci. Paris, 318 (1994),
257-260.



\bibitem{FWW} S. Fang, F.Y. Wang and B. Wu, Transportation-cost inequality on path spaces with uniform distance, {\it Stochastic Process. Appl.} 118 (2008), no. 12,
2181¨C2197.

\bibitem{FWU} S. Z. Fang,  B. Wu, \emph{Remarks on spectral gaps on
the Riemannian path  space,}  {\it arXiv:1508.07657}.


\bibitem{HN} R. Haslhofer, A. Naber, \emph{Ricci curvature and Bochner formulas for martingales,} {\it arXiv:1608.04371}.

\bibitem{Hsu0} E. P. Hsu, \emph{Logarithmic Sobolev inequalities on path
spaces over Riemannian manifolds,} Comm. Math. Phys.  189(1997),
9--16.

\bibitem{H1} E. P. Hsu, \emph{Multiplicative functional for the heat equation on manifolds with boundary,} Mich. Math. J. 50(2002),351--367.



\bibitem{N} A. Naber, Characterizations of bounded Ricci curvature on smooth and nonsmooth spaces, {\it arXiv: 1306.6512v4}.

\bibitem{TW} A. Thalmaier and F.-Y. Wang, \emph{Gradient estimates for
harmonic functions on regular domains in Riemannian manifolds,} J. Funct. Anal. 155:1(1998),109--124.

\bibitem{RS} M. R\"{o}ckner and B. Schmuland, \emph{Tightness of general $C_{1,p}$ capacities on Banach
space,} J. Funct. Anal. 108(1992), 1--12.

\bibitem{S} K.- T.Sturm, \emph{Remarks about Synthetic Upper Ricci Bounds for Metric Measure Spaces,} {\it  	arXiv:1711.01707v1}.

\bibitem{W1} F.- Y. Wang, \emph{Weak Poincar\'{e} Inequalities on path
spaces,} Int. Math. Res. Not. 2004(2004), 90--108.

\bibitem{W11}  F.-Y. Wang, \emph{Analysis on path spaces over Riemannian manifolds with boundary,} Comm. Math. Sci. 9(2011),1203--1212.

\bibitem{W2} F.- Y. Wang, \emph{Analysis for diffusion processes on Riemannian manifolds,} World Scientific, 2014.

\bibitem{W17}  F.- Y. Wang, \emph{Identifying constant curvature manifolds, Einstein manifolds, and Ricci parallel manifolds,}arXiv:1710.00276v2

\bibitem{WW} F.- Y. Wang and B. Wu, \emph{Pointwise  Characterizations  of   Curvature and  Second Fundamental Form on Riemannian Manifolds,} {\it arXiv:1605.02447}.

\bibitem{WW2} F. -Y. Wang and B. Wu, \emph{Quasi-Regular Dirichlet Forms on Free Riemannian Path and Loop Spaces,} Inf. Dimen. Anal. Quantum Probab.
and Rel. Topics 2(2009) 251--267.

\end{thebibliography}
\end{document}